\documentclass[11pt,a4paper]{article}
\usepackage[utf8]{inputenc}
\usepackage[T1]{fontenc}
\usepackage[english]{babel}
\usepackage{amsmath,amssymb,amsthm,mathtools,mathrsfs}
\usepackage{geometry}
\usepackage{microtype}\microtypesetup{expansion=false}
\usepackage[hidelinks]{hyperref}
\usepackage{cite}
\usepackage{booktabs}
\usepackage{enumitem}
\usepackage{needspace}
\usepackage{xcolor}
\allowdisplaybreaks
\newtheorem{theorem}{Theorem}[section]
\newtheorem{mainresult}{Theorem}

\newtheorem{proposition}[theorem]{Proposition}
\newtheorem{lemma}[theorem]{Lemma}
\newtheorem{corollary}[theorem]{Corollary}
\newtheorem{definition}[theorem]{Definition}

\theoremstyle{definition}
\newtheorem{remark}[theorem]{Remark}

\theoremstyle{plain}

\DeclareMathOperator{\Rad}{Rad}
\DeclareMathOperator{\Nuc}{Nuc}
\DeclareMathOperator{\GL}{GL}
\DeclareMathOperator{\Span}{span}
\DeclareMathOperator{\rank}{rank}
\DeclareMathOperator{\Aut}{Aut}
\newcommand{\F}{\mathbb{F}_2}
\newcommand{\cA}{\mathcal A}
\newcommand{\cR}{\mathcal R}
\newcommand{\aO}{\alpha_O}
\newcommand{\aM}{\alpha_M}
\newcommand{\Dtwo}{D_2}
\newcommand{\Dthree}{D_3}
\newcommand{\Ct}{\widetilde C}

\hypersetup{
 pdftitle={Counterexamples to a Cubic Nondegeneracy Conjecture},
 pdfauthor={Santiago Pineda Montoya},
 pdfsubject={Ordered separation of cubic Boolean generating functions and two families of counterexamples},
 pdfkeywords={Boolean cubic forms, Morier-Genoud, Ovsienko, twisted group algebras, ordered separation, direct sums, contraction ranks}
}

\title{\textbf{Counterexamples to a Cubic Nondegeneracy Conjecture}}
\author{\textbf{Santiago Pineda Montoya}\\
\textit{Institute of Mathematics, National University of Colombia}}
\date{}

\begin{document}
\maketitle

\begin{abstract}
Morier-Genoud and Ovsienko conjectured that a cubic Boolean generating
function satisfying an ordered separation condition on its second polar
must be equivalent to one of two symmetric normal forms.
Two infinite families of counterexamples are constructed and proved
analytically; both are direct sums of small seeds, and neither is
equivalent to the proposed forms even after linear functions are added.
The first adjoins a nondegenerate quadratic block to a ternary cubic seed
and is distinguished from the proposed forms by the radical of its third
polar.  This radical determines the associative nucleus of every twisted
quasialgebra with the given generating function, so the corresponding
algebras are not isomorphic even when the grading is forgotten.
The second family consists of direct sums of a five-variable seed and
has zero third-polar radical.  A criterion for quadratic-block
presentations establishes the separation property of the seed;
commutation-graph degrees, radical dimensions, and contraction ranks
establish non-equivalence in all dimensions of the family.
\end{abstract}

\noindent\textbf{Keywords.} Boolean cubic polynomial; twisted quasialgebra;
ordered separation; alternating trilinear form; associative nucleus;
commutation graph.\\
\textbf{2020 Mathematics Subject Classification.} Primary 17A01;
Secondary 05E18, 11E04, 15A69, 16S35.

\section{Introduction}

Cubic Boolean generating functions encode commutation and associator
signs in certain nonassociative twisted algebras graded by
$V=\mathbb F_2^n$.  Morier-Genoud and Ovsienko introduced two symmetric
families $O_n$ and $M_n$ and conjectured that an ordered nondegeneracy
condition on the second combinatorial polar characterizes their
associated generating functions up to a linear change of variables
\cite[Section~10, ``Cubic polynomials'']{MGO}.  In the terminology
adopted here, \emph{ordered separation} is the following condition:
for every pair of distinct nonzero vectors $x,y\in V$, there exists
$z\in V$ for which
\[
 \Dtwo\alpha(x,z)=1,\qquad \Dtwo\alpha(y,z)=0.
\]
It requires more than the nonvanishing of each nonzero commutation row.
The order of the two prescribed values matters.
The classification of the real forms $O_{p,q}$ under graded isomorphisms
\cite{KreuschMorierGenoud} concerns the previously introduced series;
the present question concerns arbitrary cubic Boolean generating functions.

The conjecture is disproved by two infinite constructions.  Both are
direct sums of small seeds; they differ in the invariant that separates
them from the proposed forms.  The first has a large radical in its
alternating third polar.  Consequently every twisted quasialgebra having
it as a generating function has an associative nucleus whose dimension
differs from that of every twisted quasialgebra attached to the proposed
symmetric forms.  The second construction has \emph{zero} third-polar
radical in every dimension in which it is defined.  Its five-variable
seed and all its direct sums are shown to satisfy ordered separation
analytically.  The two resulting theorems are stated below.

Throughout, $\Dtwo\alpha$ and $\Dthree\alpha$ denote the combinatorial
polars defined in Section~\ref{sec:preliminaries}; the equivalence
relations $\sim_{\GL}$ and $\sim_{\mathrm{lin}}$ are defined there as
well.  Adding a linear Boolean function does not change either polar,
so non-equivalence is proved even under the coarser second relation.
Remark~\ref{rem:literal-reading} explains the convention concerning
linear representatives in the source and the stronger equivalence
statement established here.

\Needspace{13\baselineskip}
\begin{mainresult}[Symplectic-extension counterexamples]
\label{thm:infinite-family}
For each integer $m\ge1$, the normalized cubic Boolean polynomial on
$\F^{3+2m}$
\begin{equation}
 \widehat\alpha_m(x)=x_1x_2x_3+
   \sum_{j=1}^{m}x_{2j+2}x_{2j+3}
   +\sum_{i=1}^{3+2m}x_i
 \label{eq:intro-alpha}
\end{equation}
satisfies ordered separation but is equivalent to neither $\aO$ nor
 $\aM$ in dimension $n=3+2m$, under either relation
$\sim_{\GL}$ or $\sim_{\mathrm{lin}}$.  In fact,
\[
 \dim_{\F}\Rad(\Dthree\widehat\alpha_m)=2m,
\]
whereas the third polar of either proposed normal form has zero radical
in these odd dimensions.
\end{mainresult}

The proof of Theorem~\ref{thm:infinite-family} is given in Section~\ref{sec:family}.
$\widehat\alpha_m$ is the direct sum of a normalized ternary cubic and a
normalized quadratic form with nondegenerate polar, and both summands
separate.  Section~\ref{sec:nucleus} then provides an obstruction to
\emph{ungraded} algebra isomorphisms via the associative nucleus; it does
not rely on the choice of grading or of twisting function.

A positive-dimensional third-polar radical is not essential to the
failure of the conjecture.  The second theorem gives a different
mechanism.

\begin{mainresult}[Counterexamples with zero third-polar radical]
\label{thm:radzero-family}
Let
\[
 \widehat\delta(x)=x_2x_3+x_1x_2x_5+x_1x_3x_4+
                 \sum_{i=1}^{5}x_i
 \quad (x\in\F^5),
\]
and let $\Delta_k=\widehat\delta^{\oplus k}$ be the direct sum of $k$
copies on $\F^{5k}$.  For every $k\ge1$, the polynomial $\Delta_k$ is
normalized, satisfies ordered separation, and has
\[
 \Rad(\Dthree\Delta_k)=0.
\]
Nevertheless $\Delta_k$ is equivalent to neither $\aO$ nor $\aM$
under $\sim_{\GL}$ or $\sim_{\mathrm{lin}}$ in dimension $5k$.
\end{mainresult}

The proof of Theorem~\ref{thm:radzero-family} begins with an
analytic criterion for separation in a quadratic-block presentation
and applies it to the five-variable seed.  Exact commutation-graph
degrees distinguish this seed from both symmetric forms.  Direct-sum
stability, radical dimensions, and contraction ranks then establish
non-equivalence for all $k$, with no exhaustive enumeration.
Corollary~\ref{cor:ternary-blocks} records a second radical-zero family,
built only from ternary blocks, in the dimensions $3k$ with $k\ge2$.

\subsection*{Proof strategy and organization}
The argument is organized into five stages.
\begin{enumerate}[leftmargin=2.4em,label=(\arabic*)]
\item \textbf{Definitions and invariant obstructions.}
Section~\ref{sec:preliminaries} gives the precise condition of
\cite{MGO}, normalization conventions, the two equivalence relations,
the covariance of the polars, radicals, and contraction ranks, and the
relation between graded algebra isomorphisms and $\sim_{\mathrm{lin}}$.
\item \textbf{General tools.}
Section~\ref{sec:directsums} proves that direct sums preserve ordered
separation and computes the radical and contraction ranks of the common
third polar of the proposed forms.
\item \textbf{First infinite family.}
Section~\ref{sec:seed} calculates the commutation-row structure of the
three-variable seed, and Section~\ref{sec:family} proves
Theorem~\ref{thm:infinite-family}.
\item \textbf{An obstruction independent of the grading.}
Section~\ref{sec:nucleus} computes the associative nuclei of all twisted
quasialgebras with a given generating function and deduces
nonisomorphism without requiring a grading-preserving map.
\item \textbf{Second infinite family.}
Section~\ref{sec:radzero} proves the seed's separation property by
quadratic-form arguments on affine hyperplanes and obtains the needed
graph-degree invariant by counting.  Section~\ref{sec:radzero-family}
proves Theorem~\ref{thm:radzero-family}, and
Section~\ref{sec:discussion} summarizes the two mechanisms and their
scope.
\end{enumerate}

\section{Preliminaries and the exact conjecture}
\label{sec:preliminaries}

\subsection{Boolean polynomials and their polars}

Throughout, let $V=\F^n$, where $\F$ is the field with two elements, and write $V$ additively. A Boolean polynomial on $V$ is identified with its unique algebraic normal form in
\[
 \F[x_1,\ldots,x_n]/(x_1^2-x_1,\ldots,x_n^2-x_n).
\]
Thus every monomial is square-free, and the \emph{algebraic degree} is the largest size of a monomial with nonzero coefficient.

For a function $\alpha:V\to\F$ with $\alpha(0)=0$, define its second and third combinatorial polars by
\begin{align}
 \Dtwo\alpha(x,y)
 &=\alpha(x+y)+\alpha(x)+\alpha(y),
 \label{eq:D2}\\
 \Dthree\alpha(x,y,z)
 &=\alpha(x+y+z)+\alpha(x+y)+\alpha(x+z)+\alpha(y+z)\notag\\
 &\quad+\alpha(x)+\alpha(y)+\alpha(z).
 \label{eq:D3}
\end{align}
All additions in these two formulas are in $\F$.  For every $x$, normalization
at zero gives $\Dtwo\alpha(x,0)=\Dtwo\alpha(x,x)=0$.
The operators are linear in $\alpha$.  A linear function $\ell$ has
$\Dtwo\ell=\Dthree\ell=0$, while, for $i\ne j$,
\[
 \Dtwo(x_ix_j)(x,y)=x_i y_j+x_j y_i,
 \qquad \Dthree(x_ix_j)=0.
\]
Consequently the linear part affects neither polar, and the quadratic
part affects only the second polar.  Conversely, $\Dtwo\alpha=0$ means
that $\alpha$ is additive, hence $\F$-linear.  For distinct $i,j,k$,
expansion gives
\[
 \Dthree(x_ix_jx_k)(x,y,z)
 =\sum_{\sigma\in S_3}x_{i_{\sigma(1)}}y_{i_{\sigma(2)}}z_{i_{\sigma(3)}},
 \qquad(i_1,i_2,i_3)=(i,j,k).
\]
Here $S_3$ denotes the group of permutations of the three positions. Thus, for $\deg\alpha\le3$, the third polar is an alternating trilinear form. This is the degree-three instance of combinatorial polarization \cite{Ward}; compare \cite[Proposition~4.1 and Corollary~4.2]{MGO}. Here \emph{alternating} means that the form vanishes whenever two arguments are equal. Over $\F$, every alternating trilinear form is also symmetric under all permutations of its arguments.

For an alternating trilinear form $A:V^3\to\F$, its radical is
\begin{equation}
 \Rad(A)=\{z\in V:A(z,x,y)=0\text{ for all }x,y\in V\}.
 \label{eq:radical-definition}
\end{equation}
Because $A$ is symmetric, the distinguished position of $z$ is immaterial.

\subsection{Normalized twists, twisted quasialgebras, and generating functions}

Following the terminology of \cite{MGO}, a \emph{normalized twisting function} is a map
\[
 f:V\times V\longrightarrow\F
\]
satisfying
\begin{equation}
 f(0,x)=f(x,0)=0\qquad(x\in V).
 \label{eq:normalized-f}
\end{equation}
For twisting functions, normalization means precisely Eq.~\eqref{eq:normalized-f}.
It ensures that the degree-zero basis vector is a two-sided identity.

Let $K$ be a field of characteristic different from two. The corresponding \emph{twisted group algebra} in the terminology of \cite{MGO}, or more explicitly the \emph{twisted quasialgebra}, is the cochain-twisted algebra \cite{AlbuquerqueMajid} whose underlying $K$-vector space is
\[
 K_f[V]=\bigoplus_{x\in V}K u_x
\]
with multiplication
\begin{equation}
 u_xu_y=(-1)^{f(x,y)}u_{x+y}.
 \label{eq:twisted-product}
\end{equation}
No cocycle condition is assumed. The algebra is associative precisely when $f$ is a $2$-cocycle; otherwise Eq.~\eqref{eq:twisted-product} defines a generally nonassociative unital algebra of $K$-dimension $2^n$. Each homogeneous basis element has a two-sided inverse because $u_x^2=\pm u_0$. Independently of that inverse, its left and right multiplication maps are invertible because Eq.~\eqref{eq:twisted-product} makes each map a signed permutation of the homogeneous basis. This explicit signed-permutation property, rather than a general cancellation assertion, is used below.

The commutativity and associativity defects are
\begin{align}
 \beta_f(x,y)&=f(x,y)+f(y,x),
 \label{eq:beta-from-f}\\
 \phi_f(x,y,z)&=f(x,y)+f(x+y,z)+f(y,z)+f(x,y+z).
 \label{eq:phi-from-f}
\end{align}
The notation $\partial f=\phi_f$ is used for this coboundary. They satisfy
\begin{align*}
 u_xu_y&=(-1)^{\beta_f(x,y)}u_yu_x,\\
 (u_xu_y)u_z&=(-1)^{\phi_f(x,y,z)}u_x(u_yu_z).
\end{align*}
A ternary function $\phi$ is called \emph{symmetric} if
\[
 \phi(x_{\sigma(1)},x_{\sigma(2)},x_{\sigma(3)})
 =\phi(x_1,x_2,x_3)
\]
for every permutation $\sigma\in S_3$.

\begin{definition}[Generating function {\cite[Definition~3.1]{MGO}}]
A function $\alpha:V\to\F$ with $\alpha(0)=0$ is a generating function for $K_f[V]$ if
\[
 \beta_f=\Dtwo\alpha,
 \qquad
 \phi_f=\Dthree\alpha.
\]
\end{definition}

The normalization of $f$ in Eq.~\eqref{eq:normalized-f} must not be confused with the optional \emph{basis normalization} of a generating function.  Since linear functions have zero polars, a generating function is determined by $(\beta_f,\phi_f)$ only up to the addition of a linear function, and every such class has a unique representative \cite[Lemma~6.4]{MGO} satisfying
\begin{equation}
 \alpha(e_i)=1\qquad(1\le i\le n).
 \label{eq:alpha-basis-normalization}
\end{equation}
Here a normalized generating function satisfies Eq.~\eqref{eq:alpha-basis-normalization}, whereas a normalized twist satisfies Eq.~\eqref{eq:normalized-f}. Indeed, $\Dtwo(\alpha-\alpha')=0$ says that the difference is additive and hence $\F$-linear, and prescribing its values on a basis fixes it uniquely. This freedom of a generating function must not be confused with a change of the signs of squares in a real algebra.

Morier-Genoud and Ovsienko prove that a twisted quasialgebra admits a generating function exactly when $\phi_f$ is symmetric, and that a function is a generating function of some twisted quasialgebra exactly when it is a polynomial of degree at most three \cite[Theorem~2 and Proposition~4.1]{MGO}. These results supply the original algebraic setting for the conjecture; the counterexample calculations below use $\Dtwo\alpha$ and $\Dthree\alpha$ directly.

\subsection{Ordered separation and the two equivalence relations}

In the subsection ``Cubic polynomials'' of the appendix of \cite{MGO}
(Section~10), a cubic polynomial $\alpha$ is called non-degenerate when,
for all linearly independent $x,y\in V$, there exists $z$ such that
\begin{equation}
 \Dtwo\alpha(x,z)=1,
 \qquad
 \Dtwo\alpha(y,z)=0.
 \label{eq:ordered-separation}
\end{equation}
This is the last of three notions of non-degeneracy discussed there.
The authors then state, without labelling it as a conjecture, their
belief that every non-degenerate polynomial of degree three on $V$ is
equivalent to one of the two forms \eqref{eq:aO} and \eqref{eq:aM}
below, equivalence referring to the action of $\Aut(V)\cong\GL(V)$
described in the same appendix.  No restriction on $n$ and no further
hypothesis, such as simplicity, indecomposability, or the vanishing of a
third-polar radical, is imposed.  Over $\F$, two vectors are linearly
independent exactly when they are distinct and nonzero.  The condition
is not symmetric in $x$ and $y$, and the quantifier runs over all pairs,
so it must hold for $(x,y)$ and for $(y,x)$.  To distinguish it from
other notions of nondegeneracy, this property is called \emph{ordered
separation} throughout; the adjective records the quantifier and does
not change the condition.

The two proposed normal forms, introduced in \cite[Eqs.~(4.5)--(4.6)]{MGO}, are
\begin{align}
 \aO(x)&=\sum_{i<j<k}x_ix_jx_k+
          \sum_{i<j}x_ix_j+
          \sum_i x_i,
 \label{eq:aO}\\
 \aM(x)&=\sum_{i<j<k}x_ix_jx_k+
          \sum_i x_i.
 \label{eq:aM}
\end{align}
The conventional left action of $\GL(V)$ on functions is
\[
 (T\cdot\alpha)(x)=\alpha(T^{-1}x).
\]
Elsewhere in \cite{MGO}, two generating functions with the same commutation and associator data may differ by a linear function before a basis normalization is imposed. It is therefore useful to distinguish the following relations.

\begin{definition}
The notation $\alpha\sim_{\GL}\alpha'$ means that $\alpha'=T\cdot\alpha$ for some $T\in\GL(V)$. Similarly, $\alpha\sim_{\mathrm{lin}}\alpha'$ means that
\[
 \alpha'=T\cdot\alpha+\ell
\]
for some $T\in\GL(V)$ and some linear function $\ell\in V^*$.
\end{definition}

Neither relation allows a translation $x\mapsto x+c$ of the variables.
Such a translation does not preserve $\alpha(0)=0$ in general and changes
the second polar by the contraction $\Dthree\alpha(\cdot,\cdot,c)$.

\begin{remark}[Linear representatives and equivalence conventions]
\label{rem:literal-reading}
The display of cubic polynomials in the appendix of \cite{MGO}
(Section~10, ``Cubic polynomials'') records the quadratic and cubic
terms without linear summands.  Elsewhere in that work generating
functions are considered with a prescribed basis normalization, and
the freedom to add linear functions is discussed explicitly
\cite[Lemma~6.4]{MGO}.  Since
$\Dtwo(\alpha+\ell)=\Dtwo\alpha$ for every $\ell\in V^*$,
ordered separation does not distinguish these linear representatives.
The two main theorems below prove that the displayed, basis-normalized
counterexamples cannot be taken to either symmetric form even under
the coarser relation $\sim_{\mathrm{lin}}$, which permits a linear
addition after a change of coordinates.  In particular they also rule
out equivalence under the $\GL(V)$ action used in the appendix.
This argument does not depend on how the linear representatives in the
original formulation are chosen.
For instance, on $\F^3$ the polynomial $t=x_1x_2x_3$ satisfies ordered
separation (Lemma~\ref{lem:seed-rows}) and is not $\GL(3,2)$-equivalent
to either displayed normalized form, since their numbers of nonzero
values are respectively $1$, $7$ for $\aO$, and $3$ for $\aM$.
Yet $t\sim_{\mathrm{lin}}\aM$, because
$\aM=t+x_1+x_2+x_3$ in this dimension.  Thus this small example
reflects the convention for linear representatives, whereas the main
theorems establish non-equivalence under both relations.
\end{remark}

For a trilinear form $A$ and a vector $x$, the \emph{contraction at $x$}
is the bilinear form $(y,z)\mapsto A(x,y,z)$; its rank is the rank of
its matrix in any basis.  A change of basis multiplies that matrix by
invertible matrices, so its rank is well-defined.

\begin{lemma}[Covariance and the relevant obstructions]
\label{lem:covariance}
Let $\alpha'=T\cdot\alpha+\ell$, with $T\in\GL(V)$ and
$\ell\in V^*$.  Then, for $j=2,3$,
\[
 D_j\alpha'(v_1,\ldots,v_j)
 =D_j\alpha(T^{-1}v_1,\ldots,T^{-1}v_j).
\]
In particular,
\[
 \Rad(\Dthree\alpha')=T\Rad(\Dthree\alpha).
\]
Ordered separation, the dimension of the third-polar radical, and the
multiset of ranks of contractions $\Dthree\alpha(x,\cdot,\cdot)$,
indexed by nonzero $x$, are invariant under $\sim_{\mathrm{lin}}$.
\end{lemma}

\begin{proof}
In each evaluation defining a polar, linearity gives
$T^{-1}(v_{i_1}+\cdots+v_{i_r})=
T^{-1}v_{i_1}+\cdots+T^{-1}v_{i_r}$.  The linear addition $\ell$
vanishes under both polars.  Substitution therefore proves the first
identity term by term.  A vector $z$ belongs to the radical on the left
exactly when $T^{-1}z$ belongs to the radical on the right, since
$T^{-1}x$ and $T^{-1}y$ range over all of $V$.  This gives the radical
identity and equality of its dimensions.

For contractions, let $M_x$ be the matrix of
$\Dthree\alpha(x,\cdot,\cdot)$ in a fixed basis.  The matrix for
$\alpha'$ at $Tx$ is
\[
 (T^{-1})^{\mathsf T}M_xT^{-1}.
\]
The two invertible factors preserve rank, and $x\mapsto Tx$ is a
bijection of the nonzero vectors.  Finally, a witness $z$ for a pair
$(x,y)$ becomes the witness $Tz$ for $(Tx,Ty)$.  Invertibility
preserves the condition that the pair be distinct and nonzero.  Applying
the same argument to $T^{-1}$ proves the converse implications.
\end{proof}

The finer relation $\sim_{\GL}$ is obtained by setting $\ell=0$.
Thus every obstruction in Lemma~\ref{lem:covariance} applies to both
notions of equivalence.  Equality of these invariants is not asserted
to be sufficient for equivalence.

The relation $\sim_{\mathrm{lin}}$ is also the one detected by graded
isomorphisms of twisted quasialgebras.

\begin{lemma}[Graded isomorphisms]
\label{lem:graded-iso}
Let $K$ be a field of characteristic different from two, and let $f,f'$
be normalized twisting functions on $V$ with generating functions
$\alpha$ and $\alpha'$.  Write $u_x$ and $u'_x$ for the homogeneous bases
of $K_f[V]$ and $K_{f'}[V]$.  Suppose that
$\Phi:K_f[V]\to K_{f'}[V]$ is an isomorphism of $K$-algebras and that
$T:V\to V$ is a map with $\Phi(Ku_x)=Ku'_{Tx}$ for every $x\in V$.  Then
$T\in\GL(V)$ and $\alpha'=T\cdot\alpha+\ell$ for some $\ell\in V^*$.
\end{lemma}

\begin{proof}
If $Tx=Ty$, then $\Phi(u_x)$ and $\Phi(u_y)$ are proportional, so
$u_x$ and $u_y$ are proportional because $\Phi$ is injective, and
$x=y$.  Thus $T$ is injective, hence bijective.  Write
$\Phi(u_x)=\lambda_xu'_{Tx}$ with $\lambda_x\in K^\times$.  Applying
$\Phi$ to $u_xu_y=(-1)^{f(x,y)}u_{x+y}$ gives
\[
 \lambda_x\lambda_y(-1)^{f'(Tx,Ty)}u'_{Tx+Ty}
 =(-1)^{f(x,y)}\lambda_{x+y}u'_{T(x+y)}.
\]
Both sides are nonzero multiples of basis vectors, so
$T(x+y)=Tx+Ty$.  An additive map of $\F$-vector spaces is $\F$-linear;
thus $T\in\GL(V)$.  Next apply $\Phi$ to
$u_xu_y=(-1)^{\beta_f(x,y)}u_yu_x$.  Since
$u'_{Tx}u'_{Ty}=(-1)^{\beta_{f'}(Tx,Ty)}u'_{Ty}u'_{Tx}$, this yields
\[
 (-1)^{\beta_{f'}(Tx,Ty)}\lambda_x\lambda_y\,u'_{Ty}u'_{Tx}
 =(-1)^{\beta_f(x,y)}\lambda_x\lambda_y\,u'_{Ty}u'_{Tx}.
\]
The vector $\lambda_x\lambda_yu'_{Ty}u'_{Tx}$ is nonzero, and $1\ne-1$
in $K$; hence $\beta_{f'}(Tx,Ty)=\beta_f(x,y)$, that is,
$\Dtwo\alpha'(Tx,Ty)=\Dtwo\alpha(x,y)$.  The function
$\delta=\alpha'\circ T+\alpha$ therefore satisfies $\Dtwo\delta=0$ and
$\delta(0)=0$, so it is linear.  Consequently
$\alpha'=\alpha\circ T^{-1}+\delta\circ T^{-1}=T\cdot\alpha+\ell$ with
$\ell=\delta\circ T^{-1}\in V^*$.
\end{proof}

Thus, if $\alpha\not\sim_{\mathrm{lin}}\alpha'$, no twisted quasialgebra
with generating function $\alpha$ is isomorphic to one with generating
function $\alpha'$ by an isomorphism that maps homogeneous components
onto homogeneous components.  Isomorphisms without this restriction are
treated in Section~\ref{sec:nucleus}.

\section{Direct sums and the symmetric cubic form}
\label{sec:directsums}

Both families of counterexamples are direct sums.  This section
collects two general facts.  Direct sums of
nonvacuous summands preserve ordered separation, and the common third
polar of the proposed forms has an explicit radical and explicit
contraction ranks.

If $\alpha_i:V_i\to\F$ satisfy $\alpha_i(0)=0$, define their direct sum by
\begin{equation}
 (\alpha_1\oplus\alpha_2)(x_1,x_2)
 =\alpha_1(x_1)+\alpha_2(x_2).
 \label{eq:direct-sum}
\end{equation}
The polars decompose without cross terms.  For the second polar,
\begin{align*}
 &\Dtwo(\alpha_1\oplus\alpha_2)((x_1,x_2),(y_1,y_2))\\
 &\quad=\alpha_1(x_1+y_1)+\alpha_1(x_1)+\alpha_1(y_1)
       +\alpha_2(x_2+y_2)+\alpha_2(x_2)+\alpha_2(y_2)\\
 &\quad=\Dtwo\alpha_1(x_1,y_1)+\Dtwo\alpha_2(x_2,y_2).
\end{align*}
In the seven-term expression defining the third polar, the same grouping
separates every $\alpha_1$-term from every $\alpha_2$-term.  Hence, for
$j=2,3$,
\begin{equation}
 D_j(\alpha_1\oplus\alpha_2)
 =D_j\alpha_1\oplus D_j\alpha_2.
 \label{eq:direct-sum-polars}
\end{equation}
The radical of the third polar is therefore
\begin{equation}
 \Rad\!\bigl(\Dthree(\alpha_1\oplus\alpha_2)\bigr)
 =\Rad(\Dthree\alpha_1)\oplus\Rad(\Dthree\alpha_2).
 \label{eq:direct-sum-radical}
\end{equation}
Indeed, if $(a_1,a_2)$ is in the left-hand radical, testing against
$(y_1,0),(z_1,0)$ forces $\Dthree\alpha_1(a_1,y_1,z_1)=0$ for
every $y_1,z_1$.  Testing against $(0,y_2),(0,z_2)$ gives the
analogous conclusion for $a_2$.  Conversely, if both components are
radical, each of the two terms in the decomposed third polar vanishes
for arbitrary arguments.  This gives both inclusions.  In particular,
possible cancellations for mixed arguments cannot enlarge the radical,
because arguments supported in one block test each component separately.

\begin{lemma}[Direct sums preserve ordered separation]
\label{lem:direct-sum-separation}
Let $\dim V_1,\dim V_2\ge2$ and $\alpha_i(0)=0$ for $i=1,2$.
If $\alpha_i:V_i\to\F$ satisfies ordered separation for $i=1,2$, then $\alpha_1\oplus\alpha_2$ satisfies ordered separation on $V_1\oplus V_2$.
\end{lemma}

\begin{proof}
Write $\beta_i=\Dtwo\alpha_i$.  The hypotheses imply that for any
nonzero $a\in V_i$, the row
$z\mapsto\beta_i(a,z)$ is not identically zero.  Indeed,
$\dim V_i\ge2$ permits choosing $b\notin\{0,a\}$, and ordered
separation of $(a,b)$ yields a $z$ with $\beta_i(a,z)=1$.

Now fix an ordered pair of distinct nonzero vectors
$x=(x_1,x_2)$ and $y=(y_1,y_2)$ in $V_1\oplus V_2$.
Because $x\ne y$, at least one coordinate block differs.  The
hypotheses are symmetric in the two summands, so after interchanging the
names of the summands if necessary, $x_1\ne y_1$.  This relabels
blocks only; the order of $x$ and $y$ is not changed.
For $z=(z_1,z_2)$, the second polar of the sum equals
\[
 \beta(x,z)=\beta_1(x_1,z_1)+\beta_2(x_2,z_2),
 \quad
 \beta(y,z)=\beta_1(y_1,z_1)+\beta_2(y_2,z_2).
\]
Since $x_1\ne y_1$, the case $x_1=y_1=0$ does not occur, and the
remaining possibilities are the following three.

\emph{Case I. $x_1\ne0$ and $y_1\ne0$.}
Since they differ, ordered separation in $V_1$ gives $z_1$ with
$\beta_1(x_1,z_1)=1$ and $\beta_1(y_1,z_1)=0$.
Taking $z_2=0$ makes the other summand vanish, so $\beta(x,z)=1$
and $\beta(y,z)=0$.

\emph{Case II. $x_1\ne0$ and $y_1=0$.}
By the preliminary observation, choose $z_1$ with
$\beta_1(x_1,z_1)=1$.  Set $z_2=0$.  Then the first output is one,
and the second is zero because $\beta_1(0,z_1)=0$.

\emph{Case III. $x_1=0$ and $y_1\ne0$.}
Since $x=(0,x_2)$ is nonzero, $x_2\ne0$.  In this case the second
block produces the required value one.  By the preliminary observation
applied in $V_2$, choose $z_2$ with $\beta_2(x_2,z_2)=1$, and define
$\varepsilon=\beta_2(y_2,z_2)\in\F$ (possibly $y_2=0$, in which case
$\varepsilon=0$).  The first block then compensates the value
$\varepsilon$ in the second output without affecting the first output,
because $x_1=0$.  If $\varepsilon=0$, take $z_1=0$; if
$\varepsilon=1$, use the nonempty row of $y_1$ to choose $z_1$ with
$\beta_1(y_1,z_1)=1$.  In both subcases
$\beta_1(y_1,z_1)=\varepsilon$ and $\beta_1(0,z_1)=0$.
Therefore
\[
 \beta(x,(z_1,z_2))=0+1=1,
 \qquad
 \beta(y,(z_1,z_2))=\varepsilon+\varepsilon=0.
\]

The cases exhaust $x_1\ne y_1$ and prove the assertion for either
ordering of every pair.
\end{proof}

The hypothesis $\dim V_i\ge2$ cannot be omitted.  On $V_1=\F$ every
function vanishing at zero is linear, so ordered separation holds
vacuously while the row of the nonzero vector is zero; in
$\alpha_1\oplus\alpha_2$ the vector $(1,0)$ then has an empty row, and
no pair beginning with it has a witness.

The direct-sum families are distinguished from the symmetric normal
forms by elementary invariants of alternating trilinear forms.  Let
\begin{equation}
 e_{3,n}(x)=\sum_{1\le i<j<k\le n}x_ix_jx_k,
 \qquad A_n=\Dthree e_{3,n}.
 \label{eq:e3n}
\end{equation}
Since the quadratic and linear terms of \eqref{eq:aO} and \eqref{eq:aM}
have zero third polar, $A_n$ is the common third polar of $\aO$ and
$\aM$ on $\F^n$.  Write
$s(x)=\sum_i x_i$, let $\mathbf 1=(1,\ldots,1)$, use $x\cdot y=\sum_i x_iy_i$, and let
$\operatorname{wt}(x)$ denote the Hamming weight of $x$.

\begin{lemma}[Radical and contraction ranks of the symmetric cubic form]
\label{lem:symmetric-cubic-ranks}
For $n\ge3$,
\begin{equation}
 A_n(x,y,z)
 =s(x)s(y)s(z)+s(x)(y\cdot z)+s(y)(x\cdot z)+s(z)(x\cdot y).
 \label{eq:An-formula}
\end{equation}
Consequently,
\begin{equation}
 \Rad(A_n)=
 \begin{cases}
  0,&n\text{ odd},\\
  \langle\mathbf1\rangle,&n\text{ even}.
 \end{cases}
 \label{eq:An-radical}
\end{equation}
If $n$ is odd and $x\ne0$, the alternating bilinear contraction
$A_{n,x}(y,z)=A_n(x,y,z)$ has rank over $\F$
\begin{equation}
 \rank A_{n,x}=
 \begin{cases}
  2,&\operatorname{wt}(x)\text{ even},\\
  n-1,&\operatorname{wt}(x)\text{ odd}.
 \end{cases}
 \label{eq:An-contraction-ranks}
\end{equation}
For $n=3$ the two values coincide.
\end{lemma}

\begin{proof}
\emph{Step 1. Compute the universal trilinear form.}
For $i<j<k$, the third polar of $x_ix_jx_k$ equals the sum of
its six assignments to the arguments $x,y,z$.  Summing over all
unordered triples therefore yields
\[
 A_n(x,y,z)=\sum_{i,j,k\text{ pairwise distinct}}x_iy_jz_k.
\]
Over $\F$, the indicator of three pairwise distinct indices is
$1+[i=j]+[i=k]+[j=k]$, since the expression evaluates to one
when the three indices are distinct and zero if exactly two or all
three coincide.  Summing each of the four terms over all ordered
triples gives, respectively,
\[
 s(x)s(y)s(z),\quad s(z)(x\cdot y),
 \quad s(y)(x\cdot z),\quad s(x)(y\cdot z).
\]
This proves Eq.~\eqref{eq:An-formula}.

\emph{Step 2. Compute the radical, for both parities of $n$.}
Suppose $z\in\Rad(A_n)$ and put $S=s(z)$.  For two distinct
basis vectors $e_i,e_j$, Eq.~\eqref{eq:An-formula} reduces to
\[
 0=A_n(z,e_i,e_j)=S+z_i+z_j.
\]
For any three distinct indices $i,j,k$, which exist because $n\ge3$,
compare the equations
for $(i,j)$ and $(i,k)$.  Subtracting, or equivalently adding over $\F$,
gives $z_j=z_k$.  Varying the indices shows that all coordinates
are equal to some $c\in\F$, so $z=c\mathbf1$.
The same equation now reads $S=0$; but $S=nc$ in $\F$.
If $n$ is odd this forces $c=0$.  If $n$ is even either $c=0$
or $c=1$ is possible.  It remains to check that
$\mathbf1$ is actually radical for even $n$.  Substituting
$s(\mathbf1)=n=0$ and $\mathbf1\cdot v=s(v)$ into
Eq.~\eqref{eq:An-formula} leaves two identical terms,
which cancel.  This proves Eq.~\eqref{eq:An-radical}.  For $n=2$
there are no triples, $A_2=0$, and the radical is all of $V$; this is
why $n\ge3$ is assumed.

\emph{Step 3. Determine every contraction rank for odd $n$.}
Fix $x\ne0$ and assume $n$ odd, so $s(\mathbf1)=1$.
The matrix $M_x$ of $(y,z)\mapsto A_n(x,y,z)$ in the
standard basis has zero diagonal and, for $i\ne j$, entry
$s(x)+x_i+x_j$.  Equivalently,
\begin{equation}
 M_x=s(x)(I+J)+x\mathbf1^{\mathsf T}
      +\mathbf1x^{\mathsf T},
 \qquad J=\mathbf1\mathbf1^{\mathsf T}.
 \label{eq:contraction-matrix}
\end{equation}
The diagonal of the displayed matrix is indeed zero because
$(I+J)_{ii}=0$ and $x_i+x_i=0$ over $\F$.

If $s(x)=0$, then $x$ and $\mathbf1$ are linearly independent because $\mathbf1$ has odd sum and $x$ has even sum.  In this case
\[
 M_x z=x\,s(z)+\mathbf1(x\cdot z).
\]
The two linear functionals $s(z)$ and $x\cdot z$ are
independent, because their representing vectors
$\mathbf1,x$ are independent under the nondegenerate
coordinate dot product.  Therefore the coefficient pair
$(s(z),x\cdot z)$ ranges over all of $\F^2$, and the image
of $M_x$ is exactly $\Span\{x,\mathbf1\}$.  Its rank is two.

If $s(x)=1$, the equation $M_xz=0$ becomes
\[
 z=(\mathbf1+x)s(z)+\mathbf1(x\cdot z),
\]
because $(I+J)z=z+\mathbf1s(z)$.
Applying the summation functional $s$ and using
$s(\mathbf1+x)=s(\mathbf1)+s(x)=0$ and
$s(\mathbf1)=1$, the identity $s(z)=x\cdot z$ follows.
Substituting into the displayed vector equation gives
$z=x\,s(z)$.  Hence $\ker M_x\subseteq\Span\{x\}$.
The reverse containment follows from alternation,
$A_n(x,x,y)=0$ for all $y$, so $\ker M_x=\Span\{x\}$.
Its dimension is one and the rank is $n-1$; this includes $x=\mathbf1$.
This establishes Eq.~\eqref{eq:An-contraction-ranks} in every case.
\end{proof}

\section{The three-variable seed}
\label{sec:seed}

The first construction separates the variables into two blocks. The first is the smallest genuinely cubic block,
\[
 t(u)=u_1u_2u_3
 \qquad (u\in U:=\F^3).
\]
Its role is to supply the nonzero third polar. The second block, introduced in Section~\ref{sec:family}, will be purely quadratic and will supply additional freedom in the second polar. The row structure of $t$ is therefore calculated first.

Direct expansion gives
\begin{align}
 \Dtwo t(u,r)={}&u_1u_2r_3+u_1u_3r_2+u_2u_3r_1\notag\\
 &+u_1r_2r_3+u_2r_1r_3+u_3r_1r_2.
 \label{eq:seed-beta}
\end{align}
The table below can be derived without eight unrelated substitutions.
Put $h=(1,1,1)$.  The monomial $t$ is the indicator of the single
point $h$.  If $u\ne0,h$, then $t(u)=0$ and
\[
 \Dtwo t(u,r)=\mathbf1_{\{r=h+u\}}+\mathbf1_{\{r=h\}}.
\]
The two points $h+u,h$ are distinct, so they are exactly the support.
If $u=h$, then
\[
 \Dtwo t(h,r)=1+\mathbf1_{\{r=0\}}+\mathbf1_{\{r=h\}},
\]
which is one precisely outside $\{0,h\}$.  For $u=0$ the row is zero.
These identities determine all eight row supports, recorded below.
\begin{center}
\begin{tabular}{c|l}
\toprule
$u$ & $S_u=\{r:\Dtwo t(u,r)=1\}$\\
\midrule
$000$ & $\varnothing$\\
$001$ & $\{110,111\}$\\
$010$ & $\{101,111\}$\\
$011$ & $\{100,111\}$\\
$100$ & $\{011,111\}$\\
$101$ & $\{010,111\}$\\
$110$ & $\{001,111\}$\\
$111$ & $\{001,010,011,100,101,110\}$\\
\bottomrule
\end{tabular}
\end{center}

\begin{lemma}[Seed-row properties]
\label{lem:seed-rows}
The eight functions $r\mapsto\Dtwo t(u,r)$, for $u\in U$, are pairwise
distinct. Every nonzero row is nonempty, and for every ordered pair of
distinct nonzero $u,v\in U$ there exists $r$ such that
\[
 \Dtwo t(u,r)=1,
 \qquad
 \Dtwo t(v,r)=0.
\]
Consequently $\widehat t=t+u_1+u_2+u_3$ is normalized and satisfies
ordered separation on $U$.
\end{lemma}

\begin{proof}
The table lists all $2^3$ possible values of $u$ and every $r$ at
which $\Dtwo t(u,r)=1$.  Each of the seven nonzero rows has nonempty
support.  For the six rows indexed by
$u\in\{001,010,011,100,101,110\}$, write $S_u=\{p_u,111\}$; the points
\[
 p_{001}=110,\quad p_{010}=101,\quad p_{011}=100,
 \quad p_{100}=011,\quad p_{101}=010,\quad p_{110}=001
\]
are distinct, since $p_u=h+u$.  The last row satisfies
$S_{111}=\{001,010,011,100,101,110\}$.
These six two-point supports are distinct from one another, from the
empty support at $000$, and from the six-point support at $111$.
Thus all eight row functions are distinct.

To verify \emph{ordered} separation, consider distinct nonzero
$u,v$.  If neither is $111$, take $r=p_u$; then $r\in S_u$ but
$r\notin S_v$.  If $u\ne111$ and $v=111$, take $r=111$; it belongs to
$S_u$ but not to $S_{111}$.  Finally, if $u=111$ and $v\ne111$, choose
one of the five points in $S_{111}\setminus\{p_v\}$, which is disjoint
from $S_v$.  These cases include both orders of every pair.  Finally,
$\Dtwo\widehat t=\Dtwo t$ and $\widehat t(e_i)=1$.
\end{proof}

For $n=3$, the function $\widehat t$ is exactly $\aM$ of
Eq.~\eqref{eq:aM}.

\section{Symplectic blocks and the first infinite family}
\label{sec:family}

The key idea is to adjoin a quadratic block. Quadratic terms alter $\Dtwo\alpha$ but contribute nothing to $\Dthree\alpha$. A nondegenerate symplectic block is especially useful because its linear functionals can realize any prescribed one-bit condition and any pair of independent two-bit conditions.

\begin{definition}
A bilinear form $B:W\times W\to\F$ is \emph{symplectic} if it is alternating, $B(a,a)=0$ for every $a$, and nondegenerate, meaning
\[
 \{a\in W:B(a,c)=0\text{ for every }c\in W\}=0.
\]
Equivalently, the map $W\to W^*$, $a\mapsto B(a,\cdot)$, is an isomorphism. A function $q:W\to\F$ is a \emph{quadratic refinement} of $B$ if $q(0)=0$ and $\Dtwo q=B$; such a function has degree at most two.
\end{definition}

For $m\ge1$, let $W_m=\F^{2m}$, with coordinates $a=(a_1,\ldots,a_{2m})$,
and let
\[
 B(a,c)=\sum_{j=1}^{m}\bigl(a_{2j-1}c_{2j}+a_{2j}c_{2j-1}\bigr)
\]
be the standard symplectic form.  Its standard hyperbolic quadratic
refinement is
\[
 q_0(a)=\sum_{j=1}^{m}a_{2j-1}a_{2j},
 \qquad
 \Dtwo q_0=B.
\]

\begin{lemma}[Symplectic solvability]
\label{lem:symplectic-solve}
Let $W$ be a finite-dimensional $\F$-vector space and let
$B:W\times W\to\F$ be a nondegenerate bilinear form. If $a\ne0$, then for each $\varepsilon\in\F$ there exists $c$ with $B(a,c)=\varepsilon$. If $a,b$ are distinct nonzero vectors, then for every $(\varepsilon,\eta)\in\F^2$ there exists $c$ with
\[
 B(a,c)=\varepsilon,
 \qquad
 B(b,c)=\eta.
\]
\end{lemma}

\begin{proof}
Nondegeneracy says that
\[
 \iota_B:W\longrightarrow W^*,\qquad
 a\longmapsto B(a,\cdot)
\]
is injective.  Since $W$ and $W^*$ have the same finite dimension,
it is an isomorphism.  For $a\ne0$, the functional
$\iota_B(a)$ is nonzero.  Any nonzero map from an
$\F$-vector space to the one-dimensional space $\F$ is onto.  Choose $c_1$ where its value is $1$, and use $c=0$ or $c_1$
according to the desired value $\varepsilon$.

For distinct nonzero $a,b$, the equality $a\ne b$ implies their
linear independence over the two-element field.  If
$\lambda\iota_B(a)+\mu\iota_B(b)=0$, injectivity of $\iota_B$ gives
$\lambda a+\mu b=0$, hence $\lambda=\mu=0$.
Thus $B(a,\cdot)$ and $B(b,\cdot)$ are independent and the linear map
\[
 \Lambda:W\longrightarrow\F^2,\qquad c\longmapsto(B(a,c),B(b,c))
\]
has rank two.  It is surjective; consequently every specified pair $(\varepsilon,\eta)$ has a solution.

An explicit coordinate construction is also available.  In any basis
$(g_1,\ldots,g_d)$ of $W$, the coordinates of $c$ are its coefficients,
and $B(a,c)=\sum_i c_iB(a,g_i)$.  Let
$\sigma_i=B(a,g_i)$ and $\rho_i=B(b,g_i)$.  The rows $\sigma,\rho$ are
the coordinate vectors of the independent functionals $\iota_B(a)$ and
$\iota_B(b)$, so some $2\times2$ minor is nonzero.  Thus there are indices
$i\ne j$ with $\sigma_i\rho_j+\sigma_j\rho_i=1$.  Set all coordinates of
$c$ to zero except
\[
 c_i=\varepsilon \rho_j+\eta \sigma_j,
 \qquad c_j=\varepsilon \rho_i+\eta \sigma_i.
\]
Then
\[
 B(a,c)=\varepsilon(\sigma_i\rho_j+\sigma_j\rho_i)=\varepsilon,
 \quad B(b,c)=\eta(\sigma_i\rho_j+\sigma_j\rho_i)=\eta.
\]
For a single condition, an index $i$ with $\sigma_i=1$ gives the
solution $c=\varepsilon g_i$.  These formulas make the solvability
statement constructive; they require no search over the vectors of $W$.
\end{proof}

\begin{corollary}[Symplectic blocks separate]
\label{cor:symplectic-separation}
Let $B$ be a nondegenerate alternating form on $W$, let $q$ be a
quadratic refinement of $B$, and let $\lambda\in W^*$.  Then $q+\lambda$
satisfies ordered separation on $W$.
\end{corollary}

\begin{proof}
$\Dtwo(q+\lambda)=B$.  For distinct nonzero $a,b\in W$,
Lemma~\ref{lem:symplectic-solve} gives $c$ with $B(a,c)=1$ and
$B(b,c)=0$.
\end{proof}

On $V_m=U\oplus W_m$, define
\begin{equation}
 \widehat\alpha_m=\widehat t\oplus\widehat q_m,
 \qquad
 \widehat q_m(a)=q_0(a)+\sum_{i=1}^{2m}a_i.
 \label{eq:normalized-family-abstract}
\end{equation}
In the coordinates $x=(u_1,u_2,u_3,a_1,\ldots,a_{2m})$,
\begin{equation}
 \widehat\alpha_m=x_1x_2x_3+
 \sum_{j=1}^{m}x_{2j+2}x_{2j+3}+
 \sum_{i=1}^{3+2m}x_i,
 \label{eq:family-coordinate}
\end{equation}
which is Eq.~\eqref{eq:intro-alpha}.

\begin{proof}[Proof of Theorem~\ref{thm:infinite-family}]
\emph{Step 1. Normalization and degree.}
The cubic coefficient of $x_1x_2x_3$ equals one, so
$\deg\widehat\alpha_m=3$.  Its constant coefficient is zero.
Since every square-free monomial of degree at least two vanishes on
a standard basis vector, the linear term gives
$\widehat\alpha_m(e_i)=1$ for all $i$.

\emph{Step 2. Ordered separation.}
Lemma~\ref{lem:seed-rows} proves that the normalized ternary
summand $\widehat t$ satisfies ordered separation on $U=\F^3$.
The other summand $\widehat q_m$ has second polar $B$, which is
nondegenerate, so Corollary~\ref{cor:symplectic-separation} proves
its ordered separation on $W_m=\F^{2m}$.  Both summands vanish at
zero, and their dimensions are $3\ge2$ and $2m\ge2$ because $m\ge1$.
Thus all hypotheses of Lemma~\ref{lem:direct-sum-separation} hold.
Applying that lemma to the decomposition in
Eq.~\eqref{eq:normalized-family-abstract} proves ordered separation
for $\widehat\alpha_m$.  In particular, the lemma covers every
ordered pair of distinct nonzero vectors, including vectors with
zero components in either block.

\emph{Step 3. Compute the third polar and its radical.}
Quadratic and linear polynomials have zero third polar, so
$\Dthree\widehat q_m=0$ and $\Rad(\Dthree\widehat q_m)=W_m$, while
\[
 \Dthree\widehat t(u,v,w)=\Dthree(u_1u_2u_3)(u,v,w)=\det[u\ v\ w].
\]
The middle expression expands to
\[
 u_1v_2w_3+u_1v_3w_2+u_2v_1w_3
 +u_2v_3w_1+u_3v_1w_2+u_3v_2w_1.
\]
It is the usual $3\times3$ determinant over $\F$, where the signs
of all permutations are equal to one.  If $u$ lies in
the radical of this form, then $\det[u\ v\ w]=0$ for every $v,w\in U$.
For $u\ne0$, extend $u$ to a basis $(u,v,w)$ of the three-dimensional
space $U$; its determinant is one, a contradiction.  Hence
$\Rad(\Dthree\widehat t)=0$.  By Eq.~\eqref{eq:direct-sum-radical},
\begin{equation}
 \Rad(\Dthree\widehat\alpha_m)
   =0\oplus W_m,
 \qquad
 \dim\Rad(\Dthree\widehat\alpha_m)=2m.
 \label{eq:family-radical}
\end{equation}

\emph{Step 4. Compare with the proposed forms.}
Both $\aO$ and $\aM$ have third polar $A_n$, with $n=3+2m$ odd.  By
Lemma~\ref{lem:symmetric-cubic-ranks}, $\Rad(A_n)=0$.  By
Lemma~\ref{lem:covariance}, the dimension of the third-polar radical is
invariant under $\sim_{\mathrm{lin}}$, and in particular under its finer
relation $\sim_{\GL}$.  The radical dimensions are $2m>0$ for the
constructed family and zero for both normal forms.  Neither equivalence
is possible, which completes the proof of
Theorem~\ref{thm:infinite-family}.
\end{proof}

For $m=0$ the construction gives $\widehat\alpha_0=\widehat t=\aM$ on
$\F^3$; this is why the theorem starts at $m=1$.

\begin{remark}[Other quadratic refinements]
The separation proof only uses $\Dtwo q=B$, so any quadratic refinement of any nondegenerate symplectic form on $W_m$ works, and the radical computation is unchanged. For an arbitrary refinement $q$, the basis normalization must be restored with the $q$-dependent linear function $\ell_q$ on $W_m$ defined by
\[
 \ell_q(g)=1+q(g)
\]
for every standard basis vector $g$ of $W_m$.  Equation~\eqref{eq:family-coordinate} is the canonical coordinate representative obtained from $q_0$, for which $\ell_{q_0}=\sum_ia_i$.
\end{remark}

\section{Intrinsic algebraic obstruction through the nucleus}
\label{sec:nucleus}

Theorem~\ref{thm:infinite-family} establishes non-equivalence under
coordinate changes and linear additions, and Lemma~\ref{lem:graded-iso}
converts it into the absence of graded isomorphisms.  A stronger
conclusion concerns isomorphisms of the associated nonassociative
algebras without any requirement to preserve the grading.  The dimension
of the associative nucleus provides an intrinsic obstruction to such
isomorphisms.

For a $K$-algebra $\cR$, write
\[
 [a,b,c]=(ab)c-a(bc)
\]
for the algebraic associator. The left, middle, and right nuclei consist of elements that associate with every pair in the corresponding position. Their intersection is the \emph{nucleus}
\begin{equation}
 \Nuc(\cR)=\{a:[a,x,y]=[x,a,y]=[x,y,a]=0\text{ for all }x,y\in \cR\}.
 \label{eq:nucleus-definition}
\end{equation}
It is a subspace, and it is preserved by arbitrary algebra isomorphisms.

\begin{proposition}[Nucleus criterion]
\label{prop:nucleus-criterion}
Let $K$ be a field of characteristic different from two, and let $f$ be
a normalized twisting function on $V$ with a generating function $\alpha$
of degree at most three.  Put $A=\Dthree\alpha$. Then
\begin{equation}
 \Nuc(K_f[V])
 =\Span_K\{u_z:z\in\Rad(A)\}.
 \label{eq:nucleus-criterion}
\end{equation}
Consequently,
\[
 \dim_K\Nuc(K_f[V])
 =2^{\dim_{\F}\Rad(A)}.
\]
\end{proposition}

\begin{proof}
For basis vectors the two parenthesizations can be evaluated
directly from Eq.~\eqref{eq:twisted-product}.
\begin{align*}
 (u_xu_y)u_z&=
 (-1)^{f(x,y)+f(x+y,z)}u_{x+y+z},\\
 u_x(u_yu_z)&=
 (-1)^{f(y,z)+f(x,y+z)}u_{x+y+z}.
\end{align*}
The sum of the two exponents over $\F$ is
$\phi_f(x,y,z)=A(x,y,z)$, by the definition of a generating function.
Thus
\[
 [u_x,u_y,u_z]=\varepsilon(x,y,z)u_{x+y+z},
\]
where, explicitly,
\[
 \varepsilon(x,y,z)
 =(-1)^{f(y,z)+f(x,y+z)}
      \bigl((-1)^{A(x,y,z)}-1\bigr).
\]
It is zero if $A(x,y,z)=0$ and equals $2$ or $-2$ otherwise.
The exponent arithmetic takes place in $\F$, but the coefficient
$\varepsilon$ lies in $K$.  This is exactly where the hypothesis
$\operatorname{char}K\ne2$ is needed because it guarantees that the
coefficient is nonzero whenever $A(x,y,z)=1$.

The associator is trilinear over $K$, so an element $a$ lies in the
right nucleus if and only if $[u_x,u_y,a]=0$ for all basis vectors
$u_x,u_y$; the same reduction applies to the left and middle positions.
Let $a=\sum_{z\in V}c_zu_z$ and fix $x,y$.  Then
\[
 [u_x,u_y,a]
   =\sum_{z\in V}c_z\varepsilon(x,y,z)u_{x+y+z}.
\]
For fixed $x,y$, the map $z\mapsto x+y+z$ is a bijection of $V$;
thus the vectors occurring in distinct summands are distinct basis
vectors and cannot cancel.  The displayed expression vanishes for
every $x,y$ precisely when every $z$ with $c_z\ne0$ satisfies
$A(x,y,z)=0$ for every $x,y$.  This means that the right nucleus
is the span of degrees in $\Rad(A)$.

For the other two positions,
\[
 [a,u_x,u_y]=\sum_{z\in V}c_z\varepsilon(z,x,y)u_{z+x+y},\qquad
 [u_x,a,u_y]=\sum_{z\in V}c_z\varepsilon(x,z,y)u_{x+z+y},
\]
and the same bijection argument applies.  Since $\deg\alpha\le3$, the
form $A$ is alternating trilinear over $\F$ (Section~\ref{sec:preliminaries}),
hence symmetric under permutations of its arguments, so
$\varepsilon(z,x,y)$, $\varepsilon(x,z,y)$, and $\varepsilon(x,y,z)$
vanish simultaneously.  The three nuclei therefore coincide, and their
intersection is exactly the space in Eq.~\eqref{eq:nucleus-criterion}.
The basis degrees are linearly independent and $\Rad(A)$ is an
$\F$-vector space, giving $2^{\dim_{\F}\Rad(A)}$ as the $K$-dimension.
\end{proof}

The degree hypothesis is automatic since every generating function is a
polynomial of degree at most three \cite[Proposition~4.1]{MGO}.  In the
applications below it is visible directly, since the generating
functions involved are cubic polynomials plus linear terms.

No general theorem on nuclei is needed here.  For instance, the
description~\eqref{eq:nucleus-criterion} shows directly that the nucleus
is an associative subalgebra.  The space $\Rad(A)$ is closed under addition, so
products of basis vectors of degrees in $\Rad(A)$ remain in the span,
and their associators vanish because $A$ does.

For definiteness, a twisting function with a prescribed generating
function can be written explicitly.  For a degree-at-most-three
polynomial
\[
 \alpha=\sum_i a_ix_i+\sum_{i<j}a_{ij}x_ix_j+
 \sum_{i<j<k}a_{ijk}x_ix_jx_k,
\]
the monomial rule $x_i\mapsto x_iy_i$, $x_ix_j\mapsto x_iy_j$, and
$x_ix_jx_k\mapsto x_ix_jy_k+x_iy_jx_k+y_ix_jx_k$ for $i<j<k$ of
\cite[Proposition~6.6 and rule~(6.4)]{MGO} gives the normalized twisting
function
\begin{align}
 \theta_\alpha(x,y)={}&\sum_i a_ix_iy_i+
 \sum_{i<j}a_{ij}x_iy_j\notag\\
 &+\sum_{i<j<k}a_{ijk}
 (x_ix_jy_k+x_ix_ky_j+x_jx_ky_i).
 \label{eq:canonical-twist}
\end{align}
It satisfies
\[
 \theta_\alpha(0,y)=\theta_\alpha(x,0)=0,
 \qquad
 \theta_\alpha(x,x)=\alpha(x),
 \qquad
 \partial\theta_\alpha=\Dthree\alpha,
\]
and $\theta_\alpha(x,y)+\theta_\alpha(y,x)=\Dtwo\alpha(x,y)$; in
particular $\alpha$ is a generating function for $K_{\theta_\alpha}[V]$.
These identities are verified here by the following monomial
calculations.  For the
linear and quadratic rules $x_i y_i$ and $x_i y_j$, the coboundary
vanishes by bilinearity.  For one of the cubic terms,
\begin{align*}
 &\partial(x_i x_j y_k)(x,y,z)\\
 &\quad=x_i x_j y_k+(x_i+y_i)(x_j+y_j)z_k
          +y_i y_jz_k+x_i x_j(y_k+z_k)\\
 &\quad=(x_i y_j+x_j y_i)z_k.
\end{align*}
Summing the three choices of the index in the second argument gives
all six assignments of $i,j,k$ to $x,y,z$, which is the third
polar of $x_i x_j x_k$.  On the diagonal, the same three cubic
terms give $3x_i x_jx_k=x_i x_jx_k$ in $\F$.  Symmetrization gives
the six mixed terms in $\Dtwo(x_i x_jx_k)$; the pure $x$- and
pure $y$-terms cancel in the definition of the polar.  The linear and
quadratic diagonal and symmetrization identities follow immediately
from their displayed rules.  Linearity in the coefficients $a_i,a_{ij},a_{ijk}$
therefore proves all the assertions for $\theta_\alpha$.
Over a field $K$ of characteristic different from two, put
\begin{equation}
 \cA_\alpha(K)=K_{\theta_\alpha}[V]=\bigoplus_{x\in V}Km_x,
 \qquad
 m_xm_y=(-1)^{\theta_\alpha(x,y)}m_{x+y}.
 \label{eq:canonical-algebra}
\end{equation}

The original twisting functions of $O_n$ and $M_n$
\cite[Eqs.~(1.3)--(1.4)]{MGO} consist of the cubic rule above, plus
$\sum_{i\le j}x_iy_j$ for $O_n$ and with no further term for $M_n$.
Comparison with Eq.~\eqref{eq:canonical-twist} gives
$\theta_{\aO}=f_O$ and $\theta_{\aM}=f_M+\sum_ix_iy_i$.  The term
$\sum_ix_iy_i$ is bilinear and symmetric, so it has zero coboundary and
zero symmetrization; hence $f_M$ also has generating function $\aM$.
The two twists of $M_n$ differ by a signature in the sense of
\cite[Definition~3.3]{MGO}, which changes the signs of the squares of
the generators; over $\mathbb C$ the corresponding algebras are
isomorphic as graded algebras \cite[Proposition~3.4]{MGO}.  Over a field
without a square root of $-1$ no such identification is asserted, and
none is needed.  Proposition~\ref{prop:nucleus-criterion} depends only on
the generating function.

\begin{corollary}[Ungraded nonisomorphism]
\label{cor:ungraded}
Let $m\ge1$, $n=3+2m$, $V=\F^n$, and let $K$ be a fixed field of
characteristic different from two.  Let $f$ and $f'$ be normalized
twisting functions on $V$ such that $f$ has a generating function
$\sim_{\mathrm{lin}}$-equivalent to $\widehat\alpha_m$ and $f'$ has a
generating function $\sim_{\mathrm{lin}}$-equivalent to $\aO$ or to
$\aM$.  Then $K_f[V]$ and $K_{f'}[V]$ are not isomorphic as $K$-algebras,
even though both have $K$-dimension $2^n$.  In particular
$\cA_{\widehat\alpha_m}(K)$ is isomorphic neither to
$\cA_{\aO}(K)=K_{f_O}[V]$, nor to $\cA_{\aM}(K)$, nor to $K_{f_M}[V]$.
\end{corollary}

\begin{proof}
By Lemma~\ref{lem:covariance} and Eq.~\eqref{eq:family-radical}, the
third polar of every generating function of $f$ has a radical of
dimension $2m$.  Proposition~\ref{prop:nucleus-criterion} therefore
gives
\[
 \dim_K\Nuc(K_f[V])=2^{2m}.
\]
By Lemmas~\ref{lem:covariance} and~\ref{lem:symmetric-cubic-ranks}, the
third polar of every generating function of $f'$ has zero radical in
the odd dimension $n$, so
\[
 \dim_K\Nuc(K_{f'}[V])=1.
\]
Any algebra isomorphism $\Phi:\cR\to\cR'$ preserves each associator,
because $\Phi((xy)z-x(yz))=(\Phi(x)\Phi(y))\Phi(z)-\Phi(x)(\Phi(y)\Phi(z))$.
Consequently it maps $\Nuc(\cR)$ onto $\Nuc(\cR')$ and preserves its
dimension over the fixed coefficient field $K$.  Since $2^{2m}\ne1$,
no such isomorphism exists, even without requiring compatibility
with the $V$-grading.  The final assertion follows because
$\widehat\alpha_m$, $\aO$, $\aM$, and $\aM$ are generating functions of
$\theta_{\widehat\alpha_m}$, $f_O=\theta_{\aO}$, $\theta_{\aM}$, and
$f_M$, respectively.
\end{proof}

\section{A second mechanism: a radical-zero seed}
\label{sec:radzero}

The family in Theorem~\ref{thm:infinite-family} was engineered by adjoining a large radical to the third polar. The polynomial below is \emph{not} a specialization of that theorem. It is a radical-zero counterexample in dimension five and, as Section~\ref{sec:radzero-family} will show, the seed of a second infinite family. In this mechanism the decisive information lies in the quadratic layer and hence in $\Dtwo\alpha$.

\subsection{The seed and its third polar}

Define
\begin{equation}
 \widehat\delta=x_2x_3+x_1x_2x_5+x_1x_3x_4+
 \sum_{i=1}^{5}x_i.
 \label{eq:delta}
\end{equation}
Its second polar is
\begin{align}
 \Dtwo\widehat\delta(x,z)={}&x_2z_3+x_3z_2\notag\\
 &+x_1x_2z_5+x_1x_5z_2+x_2x_5z_1
 +x_1z_2z_5+x_2z_1z_5+x_5z_1z_2\notag\\
 &+x_1x_3z_4+x_1x_4z_3+x_3x_4z_1
 +x_1z_3z_4+x_3z_1z_4+x_4z_1z_3.
 \label{eq:delta-D2}
\end{align}

\begin{lemma}[Zero third-polar radical]
\label{lem:delta-radical}
The polynomial $\widehat\delta$ satisfies
\[
 \Rad(\Dthree\widehat\delta)=0.
\]
\end{lemma}

\begin{proof}
The linear and quadratic terms of $\widehat\delta$ have zero third
polar.  With $A=\Dthree\widehat\delta$, the third polar is the
alternating three-form
\[
 A=e_1^*\wedge e_2^*\wedge e_5^*
    +e_1^*\wedge e_3^*\wedge e_4^*.
\]
Here the wedge is interpreted over $\F$; in particular, a repeated
basis index makes an alternating term vanish.  Contract at
$z=(z_1,\ldots,z_5)$ in the first argument.
\begin{align*}
 \iota_zA={}&z_1(e_2^*\wedge e_5^*+e_3^*\wedge e_4^*)
 +z_2e_1^*\wedge e_5^*+z_5e_1^*\wedge e_2^*\\
 &+z_3e_1^*\wedge e_4^*+z_4e_1^*\wedge e_3^*.
\end{align*}
The wedge products $e_i^*\wedge e_j^*$ with $i<j$ form a basis of
$\Lambda^2(V^*)$.  The coefficient of $e_2^*\wedge e_5^*$ is
$z_1$; the coefficients of $e_1^*\wedge e_5^*$,
$e_1^*\wedge e_2^*$, $e_1^*\wedge e_4^*$, and
$e_1^*\wedge e_3^*$ are, respectively, $z_2,z_5,z_3,z_4$.
If $z$ belongs to $\Rad(A)$, the contraction is zero, so all five
coordinates vanish.  Conversely $z=0$ is always radical.  Hence
$\Rad(A)=0$.
\end{proof}

\subsection{Quadratic blocks and affine restrictions}
\label{subsec:quadratic-blocks}

The separation argument becomes transparent after one coordinate is
separated from the others.  Linear summands do not affect either
polar, so write
\[
 p(\varepsilon,u)=q(u)+\varepsilon r(u),\qquad
 (\varepsilon,u)\in\F\oplus W,
\]
where $W$ is a finite-dimensional $\F$-vector space.
Here $q,r:W\to\F$ have degree at most two and $q(0)=r(0)=0$;
linear terms in either function are allowed.  Set
\[
 B=\Dtwo r,\qquad C=\Dtwo q,\qquad \Ct=C+B.
\]
These are alternating bilinear forms on $W$.  The letter $C$ in this
subsection denotes a bilinear form, not a code.  The assumption that $B$
is nondegenerate will be essential.  The next lemma explains why a
quadratic equation and a nontrivial linear equation can be solved
together.

\begin{lemma}[A quadratic function on an affine hyperplane]
\label{lem:quadratic-hyperplane}
Let $W$ be a finite-dimensional $\F$-vector space with $\dim W\ge3$, and
let $r:W\to\F$ be a function of degree at most two with $r(0)=0$ whose
polar $B=\Dtwo r$ is a nondegenerate bilinear form.  Let $L\in W^*$ be
nonzero, let $\kappa\in\F$, and put
\[
 H=\{v\in W:L(v)=\kappa\}.
\]
For every $\lambda\in W^*$ and $c\in\F$, the function
$Q=r+\lambda+c$ takes both values $0$ and $1$ on $H$.
\end{lemma}

\begin{proof}
The function $Q$ need not vanish at zero.  Since $r$ has degree at
most two and $r(0)=0$, $r(v+h)=r(v)+r(h)+B(v,h)$; adding the additive
function $\lambda$ and the constant $c$ gives the four-term identity
\[
 Q(v+h)+Q(v)+Q(h)+Q(0)=B(v,h)\qquad(v,h\in W).
\]
The zero-constant convention for $\Dtwo$ is therefore not applied
directly to $Q$ when $c\ne0$.

\emph{Step 1. The affine constraint is nonempty.}
Since $L$ is nonzero, it takes the value one somewhere and is onto
$\F$.  Choose $v_0$ with $L(v_0)=\kappa$.  Then
$H=v_0+H_0$, where $H_0=\ker L$ has dimension $\dim W-1$ by
rank--nullity.

\emph{Step 2. The polar does not vanish on the direction space.}
For a subspace $E\le W$, define
\[
 E^{\perp_B}=\{w\in W:B(w,e)=0\text{ for every }e\in E\}.
\]
The map $W\to W^*$ induced by $B$ is an isomorphism.  Under this
isomorphism $E^{\perp_B}$ corresponds to the annihilator of $E$,
which has dimension $\dim W-\dim E$.  Thus $\dim H_0^{\perp_B}=1$.
If $B(a,a')=0$ for all $a,a'\in H_0$, then
$H_0\subseteq H_0^{\perp_B}$, so $\dim W-1\le1$, contradicting
$\dim W\ge3$.
Consequently there exist $a,a'\in H_0$ with $B(a,a')=1$.

\emph{Step 3. The quadratic function takes both values on $H$.}
The function $h\mapsto Q(v_0+h)$ equals $r(h)$ plus the linear
function $B(v_0,h)+\lambda(h)$ plus a constant, so it satisfies the
same four-term identity.  Applied to $a,a'$, it gives
\[
 Q(v_0+a+a')+Q(v_0+a)+Q(v_0+a')+Q(v_0)=B(a,a')=1.
\]
All four arguments belong to $H$ because
$a,a'\in\ker L$.  Four equal values would sum to zero in $\F$,
whereas the displayed sum is one.  Hence both zero and one occur among
these four values.
This proves the assertion without enumerating the points of $W$ or $H$.
\end{proof}

A nondegenerate alternating form exists only in even dimension, so the
hypothesis means $\dim W\ge4$ in all applications.  It cannot be
lowered to $\dim W=2$.  The polar of $r(u_1,u_2)=u_1u_2$ is
nondegenerate on $\F^2$, but $r$ vanishes identically on the hyperplane
$u_1=0$.  In that case $\dim H_0=1=\dim H_0^{\perp_B}$, and Step~2 fails.

The lemma yields a criterion for a whole class of polynomials.  Its
necessity explains why some nearby cubic functions fail separation,
rather than merely providing a sufficient construction for the
particular seed.

\begin{proposition}[A quadratic-block criterion]
\label{prop:quadratic-block-criterion}
Let $W$ be a finite-dimensional $\F$-vector space with $\dim W\ge4$,
let $q,r:W\to\F$ have degree at most two with $q(0)=r(0)=0$, and assume
that $B=\Dtwo r$ is nondegenerate.  Put $C=\Dtwo q$ and $\Ct=C+B$.  Then
\[
 p(\varepsilon,u)=q(u)+\varepsilon r(u)
\]
satisfies ordered separation on $\F\oplus W$ if and only if at least
one of the alternating forms $C$ and $\Ct$ is nondegenerate.  Adding
any linear function on $\F\oplus W$ does not change this conclusion.
\end{proposition}

\begin{proof}
\emph{Step 1. Expand the row formula.}
Because $q$ and $r$ have degree at most two and vanish at zero,
$q(u+v)+q(u)+q(v)=C(u,v)$ and $r(u+v)=r(u)+r(v)+B(u,v)$; these
identities include any linear terms of $q$ and $r$.
For $x=(\varepsilon,u)$ and $z=(t,v)$,
\begin{align*}
 \Dtwo p(x,z)
 &=q(u+v)+q(u)+q(v)\\
 &\quad +(\varepsilon+t)r(u+v)+\varepsilon r(u)+t r(v)\\
 &=C(u,v)+(\varepsilon+t)
       \bigl(r(u)+r(v)+B(u,v)\bigr)
       +\varepsilon r(u)+t r(v)\\
 &=\varepsilon r(v)+C(u,v)+(\varepsilon+t)B(u,v)+t r(u).
\end{align*}
Thus
\begin{equation}
 \Dtwo p\bigl((\varepsilon,u),(t,v)\bigr)
 =\varepsilon r(v)+C(u,v)+(\varepsilon+t)B(u,v)+t r(u).
 \label{eq:block-row-formula}
\end{equation}
The cancellation occurs in $\F$.  The coefficient of $r(u)$ is $t$,
and that of $r(v)$ is $\varepsilon$.  A linear term $\mu$ of $q$ adds
the linear function $\mu(u)$ to $p$.  A linear term $\lambda$ of $r$
adds $\varepsilon\lambda(u)$, which is quadratic on $\F\oplus W$; it is
not removed by polarization, and it enters
Eq.~\eqref{eq:block-row-formula} through the values $r(u)$ and $r(v)$.

\emph{Step 2. Prove sufficiency.}
Assume that
\[
 D=C+\tau B
\]
is nondegenerate for some $\tau\in\{0,1\}$.  Fix distinct nonzero
vectors $x=(\varepsilon,u)$ and $y=(\eta,w)$.  A witness is a pair
$(t,v)$ for which the two rows in Eq.~\eqref{eq:block-row-formula}
equal $(1,0)$ in that order.  There are four cases.

\emph{Case I. $\varepsilon=\eta=0$.}
Both $u,w$ are nonzero and they are distinct, hence independent over
$\F$.  Choose $t=\tau$.  The two equations are
\[
 D(u,v)=1+\tau r(u),\qquad D(w,v)=\tau r(w).
\]
Lemma~\ref{lem:symplectic-solve}, applied to the nondegenerate form $D$,
provides $v$.  Substitution gives the first row
$1+\tau r(u)+\tau r(u)=1$ and the second row zero.

\emph{Case II. $\varepsilon=0$, $\eta=1$.}
Here $u\ne0$, while $w$ is allowed to be zero.  Choose $t=\tau$
and impose
\[
 H=\{v:D(u,v)=1+\tau r(u)\},
\]
a nonempty affine hyperplane because $D(u,\cdot)$ is a nonzero
functional.  The first row equals one on all of $H$.  The second row is
\[
 Q_y(v)=r(v)+(D+B)(w,v)+\tau r(w),
\]
that is, $r$ plus a linear function and a constant; its four-term polar
is $B$.  By Lemma~\ref{lem:quadratic-hyperplane}, $Q_y$ takes the value
zero on $H$.  Such a point gives the required outputs $(1,0)$.

\emph{Case III. $\varepsilon=1$, $\eta=0$.}
Now $w\ne0$, while $u$ is allowed to be zero.  Choose $t=\tau$ and
impose
\[
 H=\{v:D(w,v)=\tau r(w)\}.
\]
The second row vanishes on $H$, and the first row is
$r(v)+(D+B)(u,v)+\tau r(u)$.  The same hyperplane lemma supplies a
point of $H$ where this first row is one.

\emph{Case IV. $\varepsilon=\eta=1$.}
Distinctness means $u+w\ne0$; either of $u,w$ may be zero.  This time
choose $t=\tau+1$.  The two rows are
\[
 Q_u(v)=r(v)+D(u,v)+(\tau+1)r(u),\qquad
 Q_w(v)=r(v)+D(w,v)+(\tau+1)r(w).
\]
Their sum is
\[
 Q_u(v)+Q_w(v)=D(u+w,v)+(\tau+1)(r(u)+r(w)).
\]
As $D(u+w,\cdot)$ is nonzero,
\[
 H=\{v:D(u+w,v)=1+(\tau+1)(r(u)+r(w))\}
\]
is a nonempty affine hyperplane, and $Q_u+Q_w=1$ on $H$.
Lemma~\ref{lem:quadratic-hyperplane} supplies a point of $H$ where
$Q_u=1$; at that point $Q_w=0$.  This finishes all four cases and
proves sufficiency.

\emph{Step 3. Prove necessity.}
Assume instead that both $C$ and $\Ct$ are degenerate, and choose
nonzero vectors $e\in E=\Rad(C)$ and $f\in\widetilde E=\Rad(\Ct)$.  A
vector in $E\cap\widetilde E$ annihilates $B=C+\Ct$, so it is zero;
in particular $e\ne f$ and $e+f\ne0$.  For $x=(0,a)$,
Eq.~\eqref{eq:block-row-formula} gives the row
\[
 \Dtwo p\bigl((0,a),(t,v)\bigr)=
 \begin{cases}
  C(a,v),&t=0,\\
  \Ct(a,v)+r(a),&t=1.
 \end{cases}
\]
If $r(f)=1$, consider the distinct nonzero vectors $x=(0,e)$ and
$y=(0,f)$.  For $t=0$ the first row is $C(e,v)=0$; for $t=1$ the second
row is $\Ct(f,v)+r(f)=1$.  If $r(f)=0$, consider the distinct nonzero
vectors $x=(0,f)$ and $y=(0,e+f)$.  For $t=0$ the two rows are
$C(f,v)$ and $C(e+f,v)=C(f,v)$, which coincide; for $t=1$ the first
row is $\Ct(f,v)+r(f)=0$.  In either case no $(t,v)$ produces the
ordered outputs $(1,0)$, so ordered separation fails.

Finally, adding a linear function to $p$ leaves $\Dtwo p$ unchanged,
so it leaves the criterion unchanged.
\end{proof}

Only the sufficiency half is used in Theorem~\ref{thm:radzero-family}.
The necessity half shows, for instance, that replacing $x_2x_3$ by
$x_2x_5$ in Eq.~\eqref{eq:delta} destroys ordered separation without
changing the third polar.  In the notation below, $C$ becomes
$u_1v_4+u_4v_1$ and $\Ct$ becomes $u_2v_3+u_3v_2$, both of rank two.

\subsection{Applying the criterion to the seed}
\label{subsec:seed-application}

For Eq.~\eqref{eq:delta}, set $\varepsilon=x_1$ and
\begin{equation}
 q(u)=u_1u_2,\qquad r(u)=u_1u_4+u_2u_3,
 \qquad u=(x_2,x_3,x_4,x_5)\in W=\F^4.
 \label{eq:delta-quadratics}
\end{equation}
Then $\widehat\delta(\varepsilon,u)$ differs from
$q(u)+\varepsilon r(u)$ only by its linear normalization.
The relevant bilinear forms are
\begin{align}
 B(u,v)&=u_1v_4+u_4v_1+u_2v_3+u_3v_2,\notag\\
 C(u,v)&=u_1v_2+u_2v_1,\notag\\
 \Ct(u,v)&=(u_2+u_4)v_1+(u_1+u_3)v_2+u_2v_3+u_1v_4.
 \label{eq:delta-bilinear-forms}
\end{align}
In the displayed order of coordinates their matrices are
\[
 [B]=\begin{pmatrix}0&0&0&1\\0&0&1&0\\0&1&0&0\\1&0&0&0\end{pmatrix},
 \qquad
 [\Ct]=\begin{pmatrix}0&1&0&1\\1&0&1&0\\0&1&0&0\\1&0&0&0\end{pmatrix}.
\]
If $B(u,v)=0$ for every $v$, the coefficients of
$v_4,v_3,v_2,v_1$ force, successively, $u_1=u_2=u_3=u_4=0$.
If $\Ct(u,v)=0$ for every $v$, the coefficients of $v_4,v_3$
first give $u_1=u_2=0$, and those of $v_2,v_1$ then give
$u_3=u_4=0$.  Thus both forms are nondegenerate, whereas $C$ has rank
two.

\begin{proposition}[Analytic separation of the five-variable seed]
\label{prop:delta-separation}
The polynomial $\widehat\delta$ satisfies ordered separation.
\end{proposition}

\begin{proof}
The preceding coordinate calculation verifies all hypotheses of
Proposition~\ref{prop:quadratic-block-criterion}, with $\tau=1$ and
$D=\Ct$.  For clarity, the actual choices of the first witness coordinate
are $t=1$ when $(\varepsilon,\eta)\ne(1,1)$, and $t=0$ when
$(\varepsilon,\eta)=(1,1)$.  Substituting $\tau=1$ and $D+B=C$ into the
four cases of that proof gives the remaining equations
\[
\begin{array}{c|l}
 (\varepsilon,\eta)&\text{conditions on }v\\\hline
 (0,0)&\Ct(u,v)=1+r(u),\quad \Ct(w,v)=r(w),\\
 (0,1)&\Ct(u,v)=1+r(u),\quad r(v)+C(w,v)+r(w)=0,\\
 (1,0)&\Ct(w,v)=r(w),\quad r(v)+C(u,v)+r(u)=1,\\
 (1,1)&\Ct(u+w,v)=1,\quad r(v)+\Ct(u,v)=1.
\end{array}
\]
The first row of the table is solvable by independence of two linear
functionals.  In each of the other three rows, the first condition
specifies a nonempty affine hyperplane and the second is $r$ plus a linear
function and a constant, so
Lemma~\ref{lem:quadratic-hyperplane} provides the required value.
For the last row, the two original outputs differ by $\Ct(u+w,v)=1$;
therefore setting the first output to one sets the second to zero.
The linear normalization of $\widehat\delta$ changes none of these
calculations.  This completes the analytic proof of the proposition.
\end{proof}

For later degree counts, the same row formula gives the two sections
\begin{equation}
\begin{array}{c|cc}
 &t=0&t=1\\\hline
 \varepsilon=0&C(u,v)&\Ct(u,v)+r(u)\\
 \varepsilon=1&r(v)+\Ct(u,v)&r(v)+C(u,v)+r(u).
\end{array}
\label{eq:delta-slices}
\end{equation}
These formulas are specializations of Eq.~\eqref{eq:block-row-formula},
direct consequences of the second-polar formula.

\subsection{Graph degrees as exact combinatorial counts}

Following \cite[Section~5.3]{MGO}, the commutation graph joins $x$ and
$y$ when $u_x$ and $u_y$ anticommute, that is, when
$\Dtwo\alpha(x,y)=1$.  There the vertex set is all of $V$; here the
vertex $0$, which is isolated, is omitted, so $G_\alpha$ has vertex set
$V\setminus\{0\}$.  A change of coordinates induces a graph isomorphism,
while addition of a linear function leaves the graph unchanged.  Its
degree at a nonzero vertex $x$ equals $|S_x|$, where
$S_x=\{z:\Dtwo\alpha(x,z)=1\}$, since $\Dtwo\alpha(x,0)=\Dtwo\alpha(x,x)=0$
excludes the vertex $0$ and loops.  To distinguish the five-variable
seed it suffices to exhibit one degree that occurs in $G_{\widehat\delta}$
and in neither proposed graph.

\begin{proposition}[A distinguishing graph degree]
\label{prop:analytic-graph-degrees}
\begin{enumerate}[label=(\alph*),leftmargin=2em]
\item In $G_{\widehat\delta}$, the twelve vertices $x=(0,u)$ with
$(u_1,u_2)\ne(0,0)$, that is, with $x_1=0$ and $(x_2,x_3)\ne(0,0)$, have
degree sixteen.
\item For $n=5$, a vertex of Hamming weight $w$ has degree
$d_{F_O}(w)$ in $G_{\aO}$ and $d_{F_M}(w)$ in $G_{\aM}$, with the values
in Table~\eqref{eq:weight-degree-table}.  In particular, neither graph
has a vertex of degree sixteen.
\end{enumerate}
\end{proposition}

\begin{proof}
\emph{(a) The seed.}
Use the first row of Eq.~\eqref{eq:delta-slices}.  If
$x=(0,u)\ne0$, the $t=1$ slice $\Ct(u,v)+r(u)$ is a nonconstant
affine-linear function of $v$, because $\Ct$ is nondegenerate, and has
eight ones among the sixteen points of $W$.  The $t=0$ slice is
$C(u,v)=u_1v_2+u_2v_1$; it has zero ones when $u_1=u_2=0$ and eight ones
otherwise.  There are $16-4=12$ vectors $u$ with $(u_1,u_2)\ne(0,0)$, and
each gives degree $8+8=16$.

\emph{(b) The two symmetric forms.}
Their values depend only on the Hamming weight.  Define
\[
 F_O(w)=\binom w3+\binom w2+w\pmod2,\qquad
 F_M(w)=\binom w3+w\pmod2,
\]
where binomial coefficients with upper index smaller than the lower one
are zero.  Fix a vertex $x$ of weight $w$.  If a vector $z$ has $i$
ones in the support of $x$ and $j$ outside that support, then
\[
 \operatorname{wt}(z)=i+j,\qquad
 \operatorname{wt}(x+z)=w-i+j.
\]
There are $\binom wi\binom{5-w}{j}$ such vectors.  Therefore the
row degree for either weight function $F$ is the explicit binomial sum
\begin{equation}
 d_F(w)=\sum_{i=0}^{w}\sum_{j=0}^{5-w}
 \binom wi\binom{5-w}{j}
 \mathbf1_{\{F(w)+F(i+j)+F(w-i+j)=1\}},
 \label{eq:weight-degree-sum}
\end{equation}
with the condition in the indicator evaluated in $\F$.  The terms
$z=0$ $(i=j=0)$ and $z=x$ $(i=w,\ j=0)$ contribute nothing, since the
condition then reads $F(w)+F(0)+F(w)=0$.
For completeness, the partial contributions at each fixed $i$ are
listed below.  Define
\[
 K_F(w,i)=\binom wi\sum_{j=0}^{5-w}\binom{5-w}{j}
 \mathbf1_{\{F(w)+F(i+j)+F(w-i+j)=1\}}.
\]
The value sequences are
$F_O(0),\ldots,F_O(5)=(0,1,1,1,0,1)$ and
$F_M(0),\ldots,F_M(5)=(0,1,0,0,0,1)$, so the small binomial sums give
\[
\begin{array}{c|l|l}
 w&(K_{F_O}(w,i))_{i=0}^{w}&(K_{F_M}(w,i))_{i=0}^{w}\\\hline
 1&(10,10)&(10,10)\\
 2&(4,16,4)&(4,0,4)\\
 3&(1,9,9,1)&(3,3,3,3)\\
 4&(0,4,0,4,0)&(0,4,0,4,0)\\
 5&(0,0,10,10,0,0)&(0,0,10,10,0,0).
\end{array}
\]
Summing each row gives the degree values displayed below.
\begin{equation}
\begin{array}{c|r|cc|rr}
 w&\binom5w&F_O(w)&F_M(w)&d_{F_O}(w)&d_{F_M}(w)\\\hline
 1&5 &1&1&20&20\\
 2&10&1&0&24& 8\\
 3&10&1&0&20&12\\
 4&5 &0&0& 8& 8\\
 5&1 &1&1&20&20.
\end{array}
\label{eq:weight-degree-table}
\end{equation}
For example, when $w=2$, the rows of $M_5$ with $i=0,2$ each
contribute $\binom31+\binom33=4$, whereas $i=1$ contributes zero,
giving $d_{F_M}(2)=8$.  In the same case for $O_5$, the $i=0,2$
contributions are again four each; $i=1$ contributes
$\binom21\sum_{j=0}^3\binom3j=16$, giving $d_{F_O}(2)=24$.
No entry of the last two columns equals sixteen.
\end{proof}

\begin{remark}[Complete degree multisets]
The remaining degrees of $G_{\widehat\delta}$ follow from the same
slices; they are recorded for comparison and are not used below.  The
three vertices $(0,u)$ with $u_1=u_2=0\ne u$ have degree $0+8=8$.  For
$x=(1,u)$, put $j(u)=(0,0,u_1,u_2)$.  Then
$B(j(u),v)=u_2v_1+u_1v_2=C(u,v)$, $r(j(u))=0$, and
$r(u+j(u))=u_1(u_4+u_2)+u_2(u_3+u_1)=r(u)$, the two occurrences of
$u_1u_2$ cancelling in $\F$; hence $B(u+j(u),v)=\Ct(u,v)$.  The
translation identity $r(v+a)=r(v)+r(a)+B(a,v)$ turns the second row of
Eq.~\eqref{eq:delta-slices} into
\[
 r(v)+\Ct(u,v)=r(v+u+j(u))+r(u),\qquad
 r(v)+C(u,v)+r(u)=r(v+j(u))+r(u).
\]
The two products in $r(v)=v_1v_4+v_2v_3$ involve disjoint coordinate
pairs, so $r$ has $3\cdot3+1\cdot1=10$ zeros and $3\cdot1+1\cdot3=6$
ones.  Thus the ten vectors $u$ with $r(u)=0$ give degree $6+6=12$, and
the six with $r(u)=1$ give $10+10=20$.  Collecting
Table~\eqref{eq:weight-degree-table} by weight classes, the three degree
multisets are
\begin{align}
 h_{\widehat\delta}&=\{8^{(3)},12^{(10)},16^{(12)},20^{(6)}\},\notag\\
 h_{M_5}&=\{8^{(15)},12^{(10)},20^{(6)}\},\notag\\
 h_{O_5}&=\{8^{(5)},20^{(16)},24^{(10)}\},
 \label{eq:graph-histograms}
\end{align}
where $d^{(a)}$ denotes degree $d$ with multiplicity $a$.
\end{remark}

\begin{corollary}[Radical-zero counterexample]
\label{cor:delta-counterexample}
The polynomial $\widehat\delta$ is equivalent to neither $\aO$ nor
$\aM$ under $\sim_{\GL}$ or $\sim_{\mathrm{lin}}$.
\end{corollary}

\begin{proof}
Proposition~\ref{prop:delta-separation} proves ordered separation, and
Lemma~\ref{lem:delta-radical} proves zero third-polar radical.
If $\aO$ or $\aM$ were of the form $T\cdot\widehat\delta+\ell$, then by
Lemma~\ref{lem:covariance} the map $x\mapsto Tx$ would be an
isomorphism from $G_{\widehat\delta}$ onto the corresponding graph,
preserving degrees.  By
Proposition~\ref{prop:analytic-graph-degrees}, the seed has vertices of
degree sixteen, whereas neither proposed graph has such a vertex.
Neither a linear coordinate change nor an additional linear output term
can identify the generating functions.  By
Lemma~\ref{lem:graded-iso}, this also rules out isomorphisms mapping
homogeneous components onto homogeneous components between twisted
quasialgebras with these generating functions.  No assertion about
arbitrary ungraded algebra isomorphisms is made for this second
example.  By Proposition~\ref{prop:nucleus-criterion}, all three
algebras have one-dimensional nucleus.
\end{proof}

\section{The radical-zero infinite family}
\label{sec:radzero-family}

The five-variable polynomial $\widehat\delta$ is not an isolated
accident.  By Section~\ref{sec:directsums}, ordered separation is
stable under direct sums of nonvacuous summands and the third-polar
radical behaves additively.  This turns the radical-zero seed from
Section~\ref{sec:radzero} into a second infinite counterexample family.

For $k\ge1$, let
\begin{equation}
 \Delta_k=\underbrace{\widehat\delta\oplus\cdots\oplus\widehat\delta}_{k\text{ copies}}
 :\F^{5k}\longrightarrow\F.
 \label{eq:Delta-k}
\end{equation}
The block basis makes $\Delta_k(e_i)=1$ for every basis vector.

Theorem~\ref{thm:radzero-family} now follows from the analytic seed proof in
Proposition~\ref{prop:delta-separation}, its graph-degree obstruction, and
the direct-sum and contraction lemmas of Section~\ref{sec:directsums}.

\begin{proof}[Proof of Theorem~\ref{thm:radzero-family}]
\emph{Step 1. Ordered separation by direct sums.}
Proposition~\ref{prop:delta-separation} establishes ordered
separation of $\widehat\delta$ analytically.  The polynomial
vanishes at zero and is defined on a five-dimensional space.
For $k=1$ this proves the assertion.  For $k\ge2$, apply
Lemma~\ref{lem:direct-sum-separation} inductively to
$\Delta_{k-1}\oplus\widehat\delta$; its summands have dimensions
$5(k-1)\ge5$ and $5$, both vanish at zero, and both separate.
The lemma therefore proves ordered separation for every $k$.

\emph{Step 2. Compute the radical.}
Every copy of $\widehat\delta$ has zero third-polar radical by
Lemma~\ref{lem:delta-radical}.  For vectors
$x=(x_1,\ldots,x_k)$, $y=(y_1,\ldots,y_k)$,
$z=(z_1,\ldots,z_k)$ with each component in $\F^5$, the
third polar is the sum of the block polars,
\[
 \Dthree\Delta_k(x,y,z)
   =\sum_{j=1}^{k}\Dthree\widehat\delta(x_j,y_j,z_j).
\]
If $x$ lies in its radical, fix an index $j$ and choose $y,z$
supported only on that $j$th block.  Then
$\Dthree\widehat\delta(x_j,y_j,z_j)=0$ for every $y_j,z_j$,
so $x_j$ belongs to the zero radical of the seed and is zero.
This holds for every $j$, hence
\begin{equation}
 \Rad(\Dthree\Delta_k)=0.
 \label{eq:Delta-radical}
\end{equation}
The normalization propagates separately because each standard basis vector
belongs to one block, receives value one from that copy of
$\widehat\delta$, and receives zero from all other blocks.
For the degree assertion, the cubic monomial $x_1x_2x_5$ in the first
block has coefficient one; no monomial from another block has this
support.  It therefore cannot cancel.  Every summand has degree at most
three, so $\Delta_k$ has degree exactly three.  Its value at zero is zero.
These facts are independent of the radical calculation.

\emph{Step 3. Distinguish the normal forms in dimension $n=5k$.}
Both $\aO$ and $\aM$ have third polar $A_n=\Dthree e_{3,n}$.
There are three cases.

If $k=1$, Corollary~\ref{cor:delta-counterexample} already
proves non-equivalence using the commutation-graph degree of
Proposition~\ref{prop:analytic-graph-degrees}.  This special
argument is needed because the next two obstructions do not separate
the five-variable normal forms from the single seed.  All three third
polars have zero radical, and for $n=5$ the rank $n-1$ in
Eq.~\eqref{eq:An-contraction-ranks} equals four.

If $k$ is even, then $n$ is even and $n\ge10$.
Lemma~\ref{lem:symmetric-cubic-ranks} gives
$\Rad(A_n)=\Span\{\mathbf1\}$ of dimension one.
Equation~\eqref{eq:Delta-radical} gives radical dimension zero
for $\Delta_k$.  More explicitly, a hypothetical identity
$\aO=T\cdot\Delta_k+\ell$ or
$\aM=T\cdot\Delta_k+\ell$ would imply
$\Rad(A_n)=T\Rad(\Dthree\Delta_k)=0$ by
Lemma~\ref{lem:covariance}, contradicting
$\Rad(A_n)=\Span\{\mathbf1\}$ of dimension one.

Finally, if $k\ge3$ is odd, then $n=5k\ge15$ is odd.  Choose
the nonzero vector $x$ that equals $e_1$ in the first five-variable
block and zero in every other block.  The contraction of
$\Dthree\Delta_k$ at $x$ is supported only in that first block.
Because the cubic part of the seed is
$x_1x_2x_5+x_1x_3x_4$, the contraction is
\[
 \iota_x(\Dthree\Delta_k)
   =e_2^*\wedge e_5^*+e_3^*\wedge e_4^*.
\]
In the ordered basis $e_2,e_5,e_3,e_4$, its nonzero block is
\[
 \begin{pmatrix}
 0&1&0&0\\1&0&0&0\\0&0&0&1\\0&0&1&0
 \end{pmatrix},
\]
whose determinant is one because it is block diagonal with two copies of
$\left(\begin{smallmatrix}0&1\\1&0\end{smallmatrix}\right)$.
The direction $e_1$ in the first block and every coordinate in the
other $k-1$ blocks are in the kernel, and the displayed four directions
are independent.  There are no off-diagonal couplings to the other
blocks because the polynomial is a direct sum.  Thus the rank of the
full $5k\times5k$ contraction matrix is exactly four, not merely at
most four.  The vector $e_2$ of the same block would not serve because its
contraction is $e_1^*\wedge e_5^*$, of rank two.  In contrast,
Eq.~\eqref{eq:An-contraction-ranks}
says every nonzero contraction of $A_n$ has rank $2$ or
$n-1$, the latter at least $14$.  Neither value equals four.

\emph{Step 4. Pass to both equivalence relations.}
For any $T\in\GL(n,2)$ and $\ell\in V^*$, the third polar of
$T\cdot\Delta_k+\ell$ is obtained by the simultaneous change of
variables $T^{-1}$ in all three entries of $\Dthree\Delta_k$.
This preserves radical dimension and the multiset of contraction
ranks.  It also preserves the commutation graph up to the relabeling
$x\mapsto Tx$, because the second polar is equivariant and
$\Dtwo\ell=0$.  Hence the obstruction appropriate to each of the
three cases rules out $\sim_{\mathrm{lin}}$, and a fortiori rules out
$\sim_{\GL}$.  This completes the proof for every $k\ge1$.
\end{proof}

The same lemmas give a second radical-zero family from the ternary seed
alone.  It does not reach dimension five, where
Theorem~\ref{thm:radzero-family} starts.

\begin{corollary}[Ternary blocks]
\label{cor:ternary-blocks}
For $k\ge2$, the normalized polynomial
\[
 \Gamma_k(x)=\sum_{j=0}^{k-1}x_{3j+1}x_{3j+2}x_{3j+3}+\sum_{i=1}^{3k}x_i
 =\widehat t^{\oplus k}
\]
on $\F^{3k}$ satisfies ordered separation and has
$\Rad(\Dthree\Gamma_k)=0$, but is equivalent to neither $\aO$ nor $\aM$
under $\sim_{\GL}$ or $\sim_{\mathrm{lin}}$.
\end{corollary}

\begin{proof}
Lemma~\ref{lem:seed-rows} and induction with
Lemma~\ref{lem:direct-sum-separation} give ordered separation, since
each summand has dimension three.  By Step~3 of the proof of
Theorem~\ref{thm:infinite-family}, $\Rad(\Dthree\widehat t)=0$, and
Eq.~\eqref{eq:direct-sum-radical} gives $\Rad(\Dthree\Gamma_k)=0$.  If
$k$ is even, then $n=3k$ is even and $\Rad(A_n)$ has dimension one by
Lemma~\ref{lem:symmetric-cubic-ranks}.  If $k\ge3$ is odd, let $x$ be
nonzero in exactly two blocks, with components $x^{(1)},x^{(2)}$.  The
contraction of $\Dthree\Gamma_k$ at $x$ is block diagonal with blocks
$(y,z)\mapsto\det[x^{(i)}\ y\ z]$ for $i=1,2$ and zero elsewhere.  Each
such block is an alternating form on $\F^3$ whose radical is
$\Span\{x^{(i)}\}$, because $\det[x^{(i)}\ y\ z]=0$ for all $z$ exactly
when $y\in\Span\{x^{(i)}\}$; hence each block has rank two and the
contraction has rank four.  Every nonzero contraction of $A_n$ has rank
$2$ or $n-1\ge8$.  Lemma~\ref{lem:covariance} concludes both cases.
For $k=1$, $\Gamma_1=\widehat t=\aM$ on $\F^3$.
\end{proof}

\section{Discussion and scope}
\label{sec:discussion}

The two constructions disprove the ordered-separation characterization
proposed in \cite{MGO}.  Both are direct sums of small seeds, and
Lemma~\ref{lem:direct-sum-separation} is the common source of ordered
separation in higher dimensions; what differs is the invariant that
separates them from the symmetric normal forms.  In the
symplectic-extension family the ternary seed carries the whole third
polar, while the quadratic block contributes a radical of dimension
$2m$ together with the linear freedom needed for ordered separation.
That radical determines the associative nucleus of every twisted
quasialgebra with this generating function
(Corollary~\ref{cor:ungraded}), an invariant under
arbitrary algebra isomorphisms over a fixed coefficient field of
characteristic different from two.

The family $\Delta_k$ has identically zero third-polar radical.  Its
seed is distinguished from the symmetric normal forms using
commutation-graph degrees, while higher direct sums are distinguished
by radical dimension or contraction rank.  Together with
Corollary~\ref{cor:ternary-blocks}, this shows that the conjecture
cannot be restored merely by imposing triviality of the third-polar
radical.

Both the seed separation criterion and the infinite-family arguments
are analytic.  The proofs require no exhaustive orbit classification
or numerical verification.  The constructions do not amount to a
classification of all cubic Boolean functions satisfying ordered
separation, and direct sums are not claimed to account for all of them.
Such a classification would need to retain information from the
quadratic part of the generating function, because the third polar alone
does not determine ordered separation; the example following
Proposition~\ref{prop:quadratic-block-criterion} illustrates this.  In
particular, the structure of the stabilizer action on quadratic
extensions of a fixed alternating three-form remains a separate
question.


\begin{thebibliography}{99}
\bibitem{MGO}
S. Morier-Genoud and V. Ovsienko,
\emph{A Series of Algebras Generalizing the Octonions and Hurwitz--Radon Identity},
Commun. Math. Phys. \textbf{306}, 83--118 (2011).
Statement and equation numbers are cited from the preprint version
arXiv:1003.0429v4.

\bibitem{KreuschMorierGenoud}
M. Kreusch and S. Morier-Genoud,
\emph{Classification of the Algebras $O_{p,q}$},
Commun. Algebra \textbf{43}, no.~9, 3799--3815 (2015).

\bibitem{Ward}
H. N. Ward,
\emph{Combinatorial Polarization},
Discrete Math. \textbf{26}, 185--197 (1979).

\bibitem{AlbuquerqueMajid}
H. Albuquerque and S. Majid,
\emph{Quasialgebra Structure of the Octonions},
J. Algebra \textbf{220}, 188--224 (1999).
\end{thebibliography}
\end{document}